\documentclass[12pt]{article}%
\usepackage{amsfonts}
\usepackage{sw20bams}
\usepackage{amsmath}
\usepackage{amssymb}
\usepackage{graphicx}%
\providecommand{\U}[1]{\protect\rule{.1in}{.1in}}
\begin{document}

\title{M/M/$c$ Queues and the Poisson Clumping Heuristic}
\author{Steven Finch}
\date{April 8, 2019}
\maketitle

\begin{abstract}
In continuous time, customers arrive at random. \ Each waits until one of $c$
servers is available; each thereafter departs at random. \ The distribution of
maximum line length of idle customers was studied over 25 years ago. \ We
revisit two good approximations of this, employing a discrete Gumbel
formulation and detailed graphics to describe simulation outcomes.

\end{abstract}

\footnotetext{Copyright \copyright \ 2019 by Steven R. Finch. All rights
reserved.}Consider an M/M/$c$ queue with arrival rate $\lambda$ and service
rate $\mu$. \ Let $M_{n}$ denote the maximum queue length over the time
interval $[0,n]$. \ For integer $k$, we could study $\mathbb{P}\left\{
M_{n}<k\right\}  $ asymptotically as a function of $n$, as was done in
\cite{Ald-heu} for the case $c=1$.\ \ We prefer, however, to suppress the
dependence on $n$\ somewhat, separating (in essence)\ signal from noise. \ Let
$k=\log_{c\mu/\lambda}(n)+h+1$ as defined in \cite{Fi1-heu}.\ \ The Poisson
clumping heuristic asserts that, if $\lambda<c\,\mu$, then%
\[
\mathbb{P}\left\{  M_{n}\leq\log_{c\mu/\lambda}(n)+h\right\}  \sim\exp\left[
-\frac{c^{c-2}\lambda\,\mu^{c-3}(c\mu-\lambda)^{2}}{\sum\limits_{j=1}%
^{c}j!\tbinom{c-1}{j-1}\lambda^{c-j}\mu^{j-1}}\left(  \frac{\lambda}{c\mu
}\right)  ^{h+1}\right]
\]
as $n\rightarrow\infty$. \ The finite sum involving binomial cofficients is
explained in Section 1. In particular,%
\[
\mathbb{P}\left\{  M_{n}\leq\log_{\mu/\lambda}(n)+h\right\}  \sim\exp\left[
-\frac{\lambda(\mu-\lambda)^{2}}{\mu^{2}}\left(  \frac{\lambda}{\mu}\right)
^{h+1}\right]
\]
for $c=1$,%
\[
\mathbb{P}\left\{  M_{n}\leq\log_{2\mu/\lambda}(n)+h\right\}  \sim\exp\left[
-\frac{\lambda(2\mu-\lambda)^{2}}{\mu(2\mu+\lambda)}\left(  \frac{\lambda
}{2\mu}\right)  ^{h+1}\right]
\]
for $c=2$,%
\[
\mathbb{P}\left\{  M_{n}\leq\log_{3\mu/\lambda}(n)+h\right\}  \sim\exp\left[
-\frac{3\lambda(3\mu-\lambda)^{2}}{6\mu^{2}+4\lambda\mu+\lambda^{2}}\left(
\frac{\lambda}{3\mu}\right)  ^{h+1}\right]
\]
for $c=3$,%
\[
\mathbb{P}\left\{  M_{n}\leq\log_{4\mu/\lambda}(n)+h\right\}  \sim\exp\left[
-\frac{16\lambda\mu(4\mu-\lambda)^{2}}{24\mu^{3}+18\lambda\mu^{2}+6\lambda
^{2}\mu+\lambda^{3}}\left(  \frac{\lambda}{4\mu}\right)  ^{h+1}\right]
\]
for $c=4$ and%
\[
\mathbb{P}\left\{  M_{n}\leq\log_{5\mu/\lambda}(n)+h\right\}  \sim\exp\left[
-\frac{125\lambda\mu^{2}(5\mu-\lambda)^{2}}{120\mu^{4}+96\lambda\mu
^{3}+36\lambda^{2}\mu^{2}+8\lambda^{3}\mu+\lambda^{4}}\left(  \frac{\lambda
}{5\mu}\right)  ^{h+1}\right]
\]
for $c=5$. \ Also,
\begin{align*}
\mathbb{E}\left(  M_{n}\right)   &  \approx\frac{\ln(n)}{\ln(\frac{\mu
}{\lambda})}+\frac{\gamma+\ln\left(  \frac{\lambda^{2}(\mu-\lambda)^{2}}%
{\mu^{3}}\right)  }{\ln(\frac{\mu}{\lambda})}+\frac{1}{2}\\
&  \approx(2.4663034623...)\ln(n)-(7.2049448811...)
\end{align*}
for $(c,\lambda,\mu)=(1,1/3,1/2)$,%
\begin{align*}
\mathbb{E}\left(  M_{n}\right)   &  \approx\frac{\ln(n)}{\ln(\frac{2\mu
}{\lambda})}+\frac{\gamma+\ln\left(  \frac{\lambda^{2}(2\mu-\lambda)^{2}}%
{2\mu^{2}(2\mu+\lambda)}\right)  }{\ln(\frac{2\mu}{\lambda})}+\frac{1}{2}\\
&  \approx(2.4663034623...)\ln(n)-(6.7552845943...).
\end{align*}
for $(c,\lambda,\mu)=(2,1/3,1/4)$,%
\begin{align*}
\mathbb{E}\left(  M_{n}\right)   &  \approx\frac{\ln(n)}{\ln(\frac{3\mu
}{\lambda})}+\frac{\gamma+\ln\left(  \frac{\lambda^{2}(3\mu-\lambda)^{2}}%
{\mu\left(  6\mu^{2}+4\lambda\mu+\lambda^{2}\right)  }\right)  }{\ln
(\frac{3\mu}{\lambda})}+\frac{1}{2}\\
&  \approx(2.4663034623...)\ln(n)-(6.2049448811...)
\end{align*}
for $(c,\lambda,\mu)=(3,1/3,1/6)$,%
\begin{align*}
\mathbb{E}\left(  M_{n}\right)   &  \approx\frac{\ln(n)}{\ln(\frac{4\mu
}{\lambda})}+\frac{\gamma+\ln\left(  \frac{4\lambda^{2}(4\mu-\lambda)^{2}%
}{24\mu^{3}+18\lambda\mu^{2}+6\lambda^{2}\mu+\lambda^{3}}\right)  }{\ln
(\frac{4\mu}{\lambda})}+\frac{1}{2}\\
&  \approx(2.4663034623...)\ln(n)-(5.6015876099...)
\end{align*}
for $(c,\lambda,\mu)=(4,1/3,1/8)$ and%
\begin{align*}
\mathbb{E}\left(  M_{n}\right)   &  \approx\frac{\ln(n)}{\ln(\frac{5\mu
}{\lambda})}+\frac{\gamma+\ln\left(  \frac{25\lambda^{2}\mu(5\mu-\lambda)^{2}%
}{120\mu^{4}+96\lambda\mu^{3}+36\lambda^{2}\mu^{2}+8\lambda^{3}\mu+\lambda
^{4}}\right)  }{\ln(\frac{5\mu}{\lambda})}+\frac{1}{2}\\
&  \approx(2.4663034623...)\ln(n)-(4.9642624490...)
\end{align*}
for $(c,\lambda,\mu)=(5,1/3,1/10)$, where $\gamma$ denotes Euler's constant
\cite{Fn-heu}. \ With regard to expected maximums, in a hospital emergency
room ($\lambda=1/3$), one fast doctor ($\mu=1/2$) outperforms $c$ very slow
doctors (each $\mu=1/(2c)$).

A\ higher-order approximation for the probability is \cite{Sfz-heu, MP-heu,
HM-heu}%
\[
\mathbb{P}\left\{  M_{n}\leq\log_{c\mu/\lambda}(n)+h\right\}  \sim\exp\left[
-\frac{n\,c^{c-2}\lambda^{h+c+1}\mu^{c-3}(c\mu-\lambda)^{2}}{\left\{
n\,\lambda^{c-1}(c\mu)^{h+1}-\lambda^{h+1}(c\mu)^{c-1}\right\}  \sum
\limits_{j=1}^{c}j!\tbinom{c-1}{j-1}\lambda^{c-j}\mu^{j-1}}\right]
\]
with the same relation between (real) $h$ and (integer) $k$ as earlier. \ The
dependence of the probability on $n$ is more visible here; allowing
$n\rightarrow\infty$ within the square brackets yields exactly the same
expression as before.

In greater detail, Serfozo \cite{Sfz-heu} examined the distribution of
$M_{\nu}$, the maximum queue length over $\ell$ busy cycles, where $\nu$ is
the $\ell^{\text{th}}$ time the system becomes empty. \ Very important
corrections to Serfozo's second table appeared in \cite{HM-heu}; observe the
existence of exact probabilistic results in this special case. \ McCormick
\&\ Park \cite{MP-heu} examined the distribution of $M_{n}$ for arbitrary $n$,
following \cite{Sfz-heu}. \ A\ missing coefficient $c!/c^{c}$ in McCormick
\&\ Park's formula (2.22) was uncovered in \cite{HM-heu}; no exact results are
generally available here.

We simulated $10^{6}$ M/M/$c$ queues for each of the following choices of
$(n,c,\lambda,\mu)$:

\begin{itemize}
\item $n=1000$ or $2500$

\item $c=1$, $2$, $3$, $4$ or $5$

\item $\lambda=1/3$

\item $\mu=1/(2c)$
\end{itemize}

\noindent and indicate both low-order approximation (red) and high-order
approximation (green) histograms against empirical outcomes (blue). \ The
green segments are always closer to the blue segments than the red segments.
\ Also, the discrepancies become smaller as $n$ grows larger. \ We used
clumping heuristic-based estimates for the mean (from earlier) and likewise
\[
\mathbb{V}\left(  M_{n}\right)  \approx\frac{\pi^{2}}{6}\frac{1}{\ln\left(
\frac{c\mu}{\lambda}\right)  ^{2}}+\frac{1}{12}%
\]
for the variance. \ The fit between moments is surprisingly good; no attempt
was made to employ a more sophisticated underlying formula.%
\begin{figure}[ptb]%
\centering
\includegraphics[
height=3.0234in,
width=3.4541in
]%
{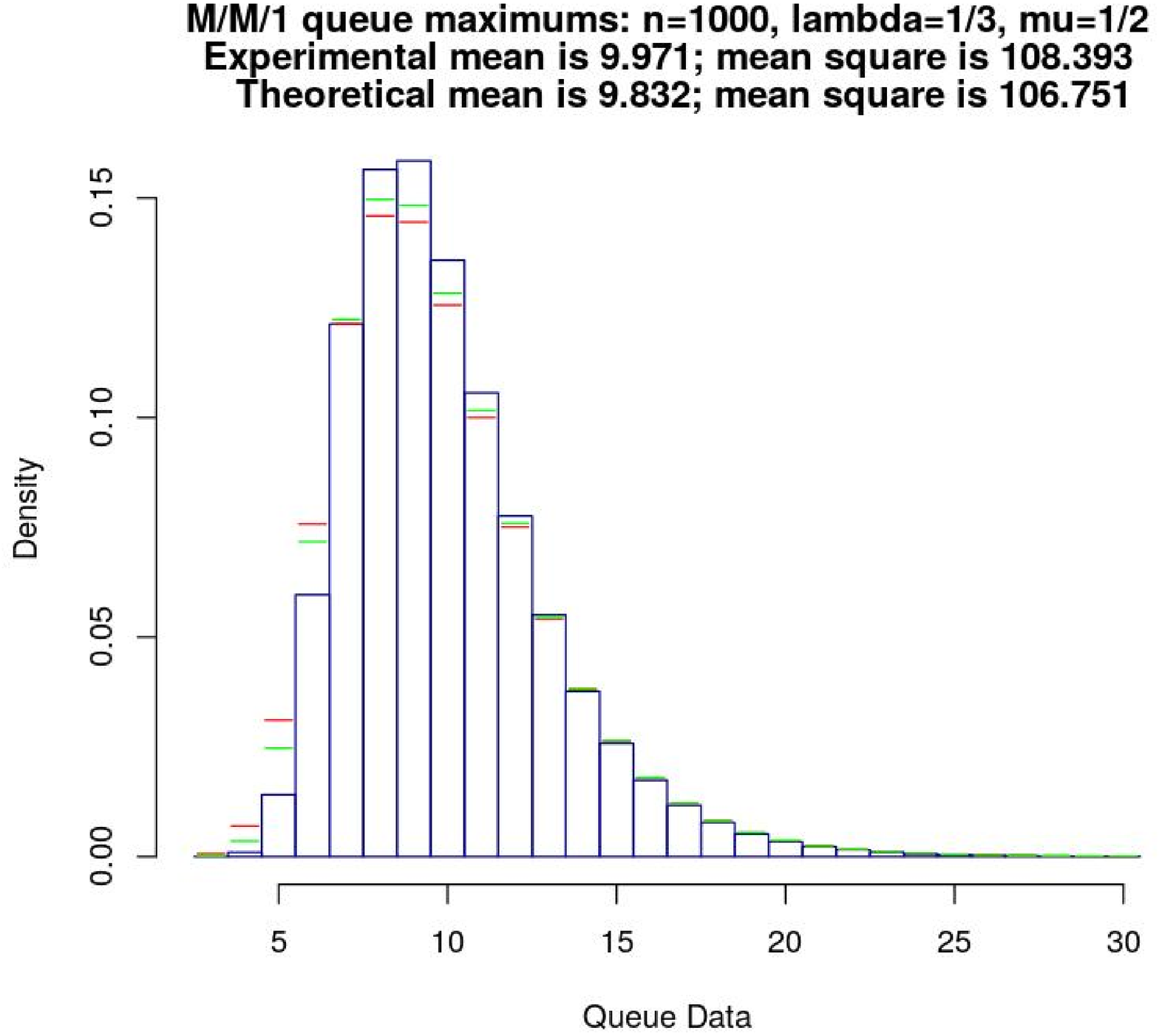}%
\end{figure}
\begin{figure}[ptb]%
\centering
\includegraphics[
height=3.0234in,
width=3.4541in
]%
{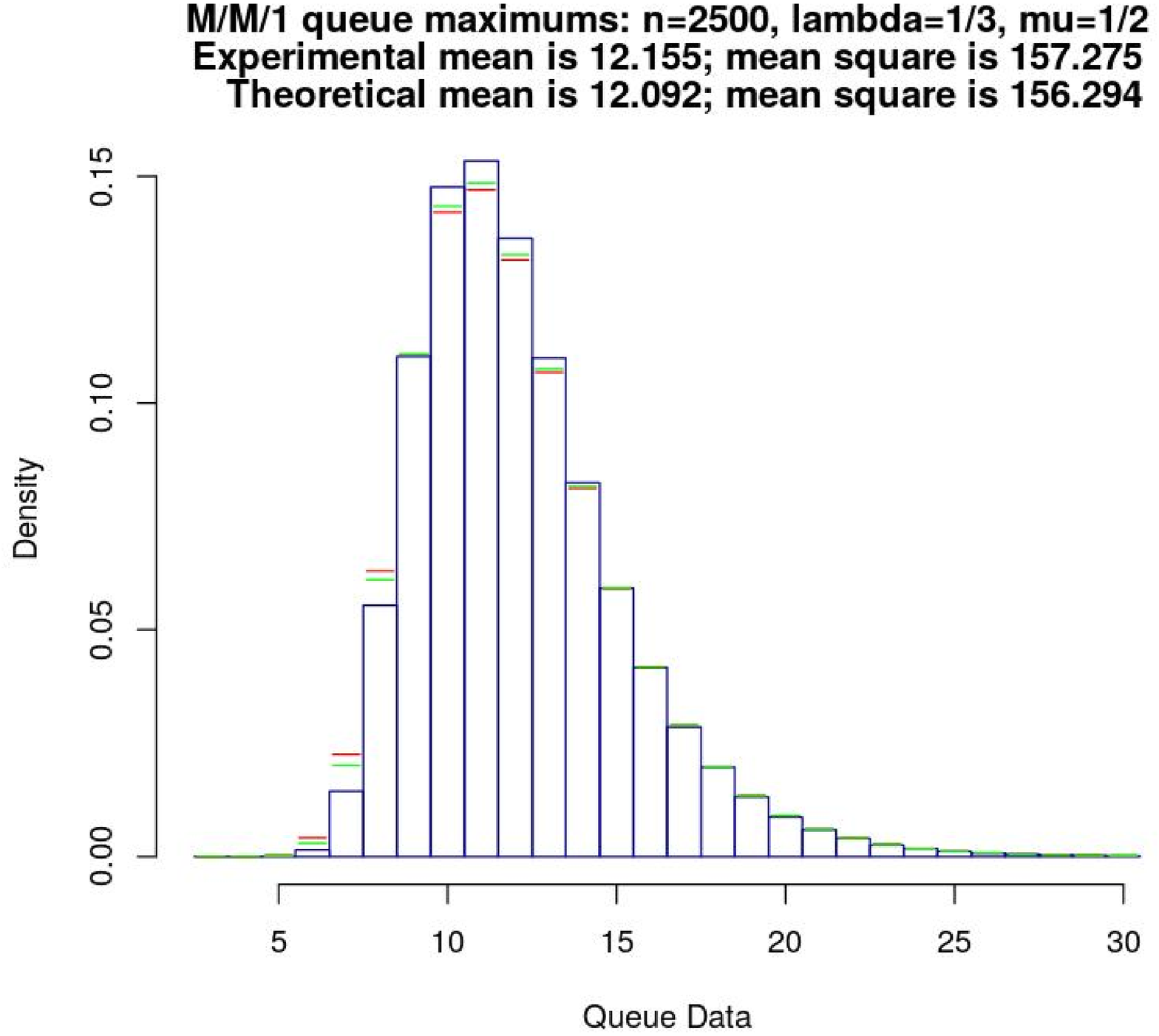}%
\end{figure}
\begin{figure}[ptb]%
\centering
\includegraphics[
height=3.0234in,
width=3.4541in
]%
{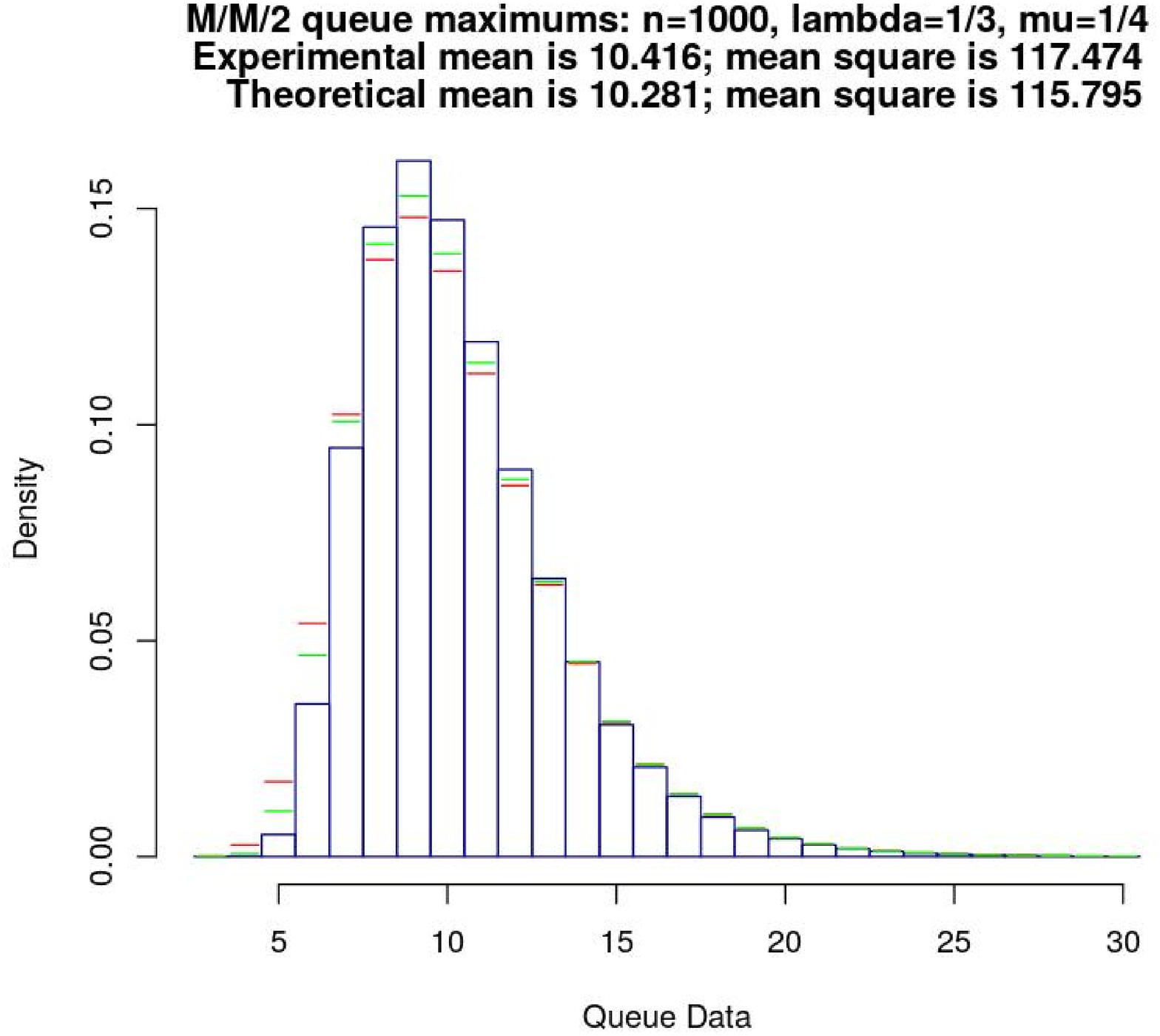}%
\end{figure}
\begin{figure}[ptb]%
\centering
\includegraphics[
height=3.0234in,
width=3.4541in
]%
{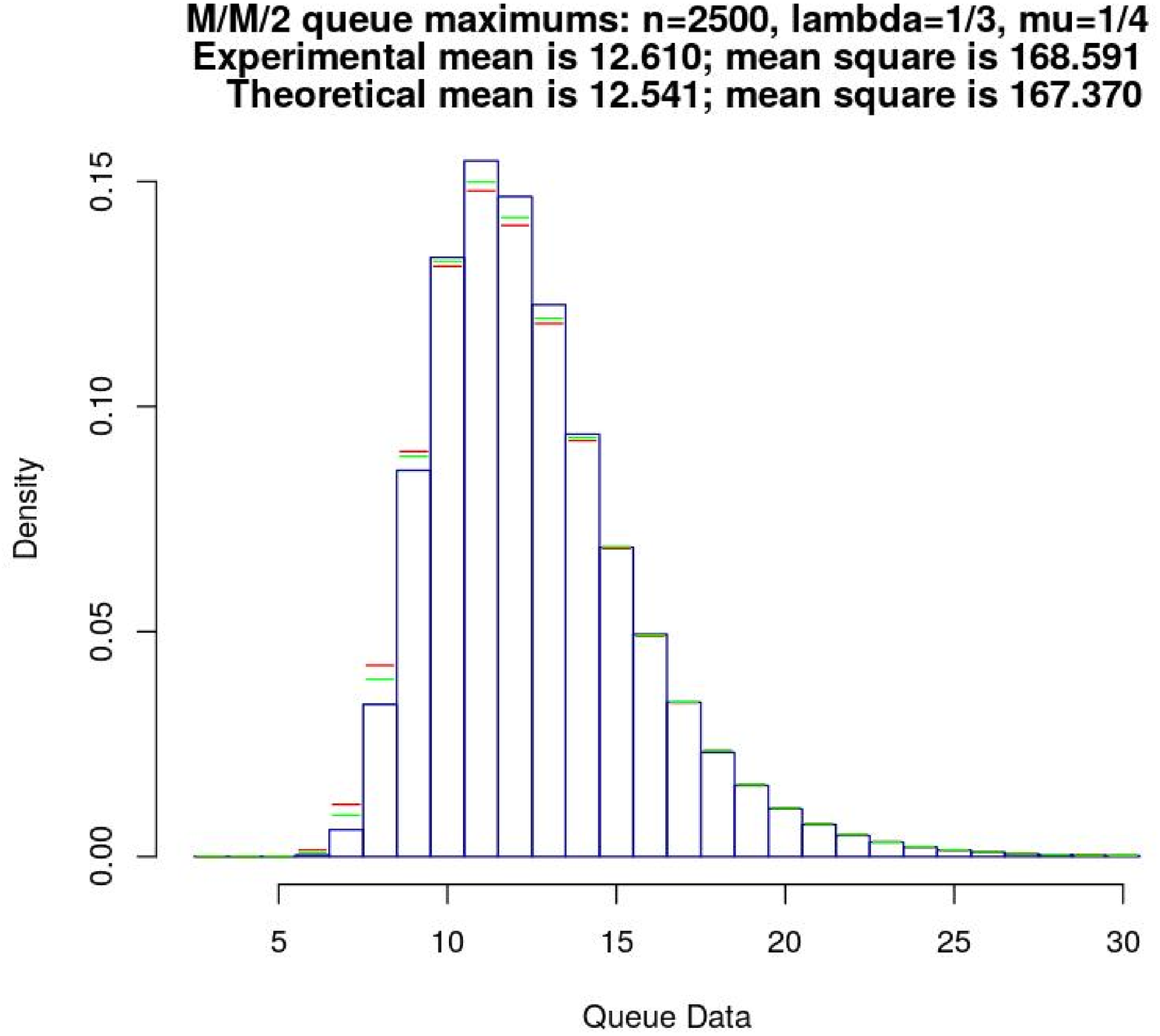}%
\end{figure}
\begin{figure}[ptb]%
\centering
\includegraphics[
height=3.0234in,
width=3.4541in
]%
{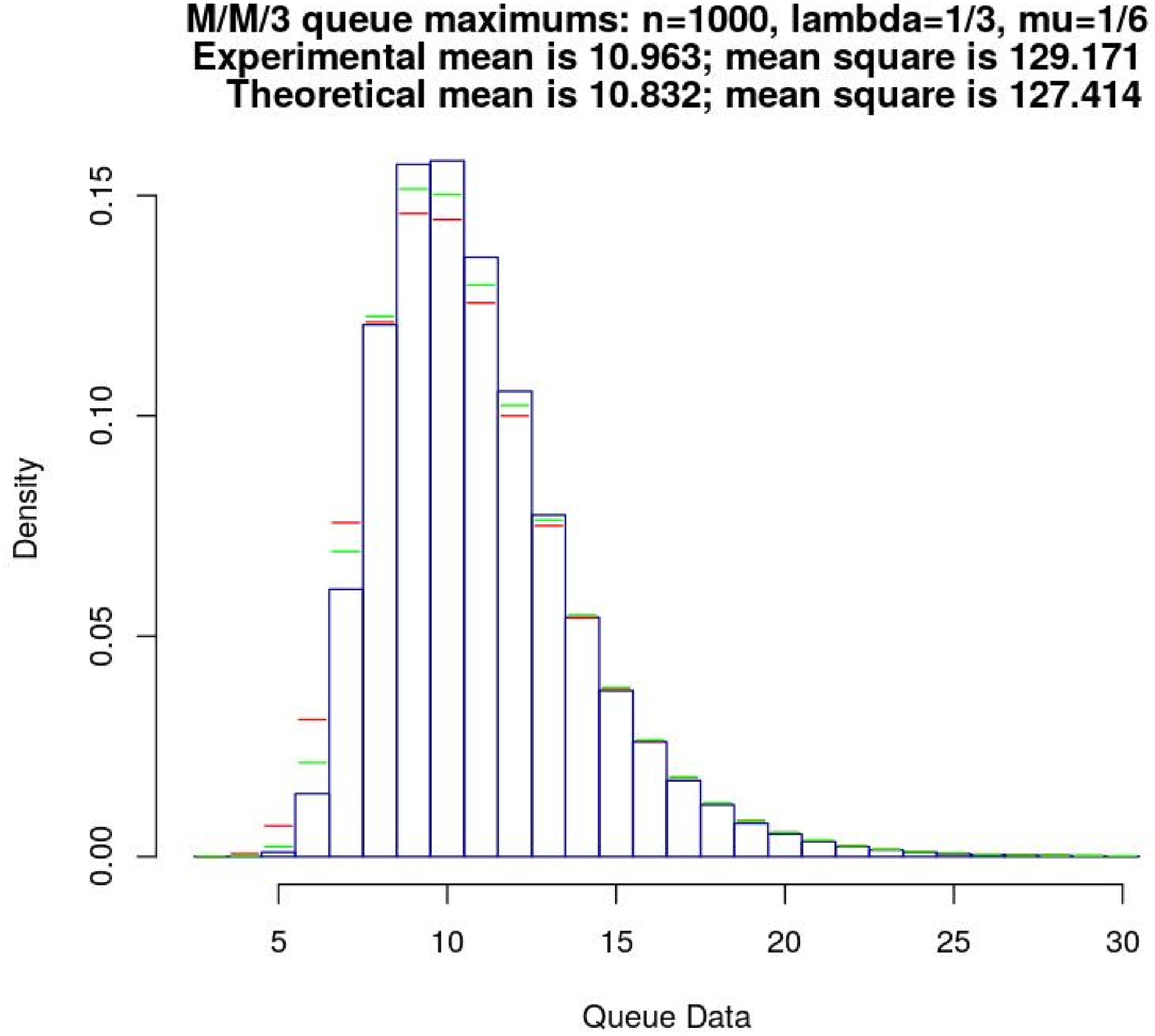}%
\end{figure}
\begin{figure}[ptb]%
\centering
\includegraphics[
height=3.0234in,
width=3.4541in
]%
{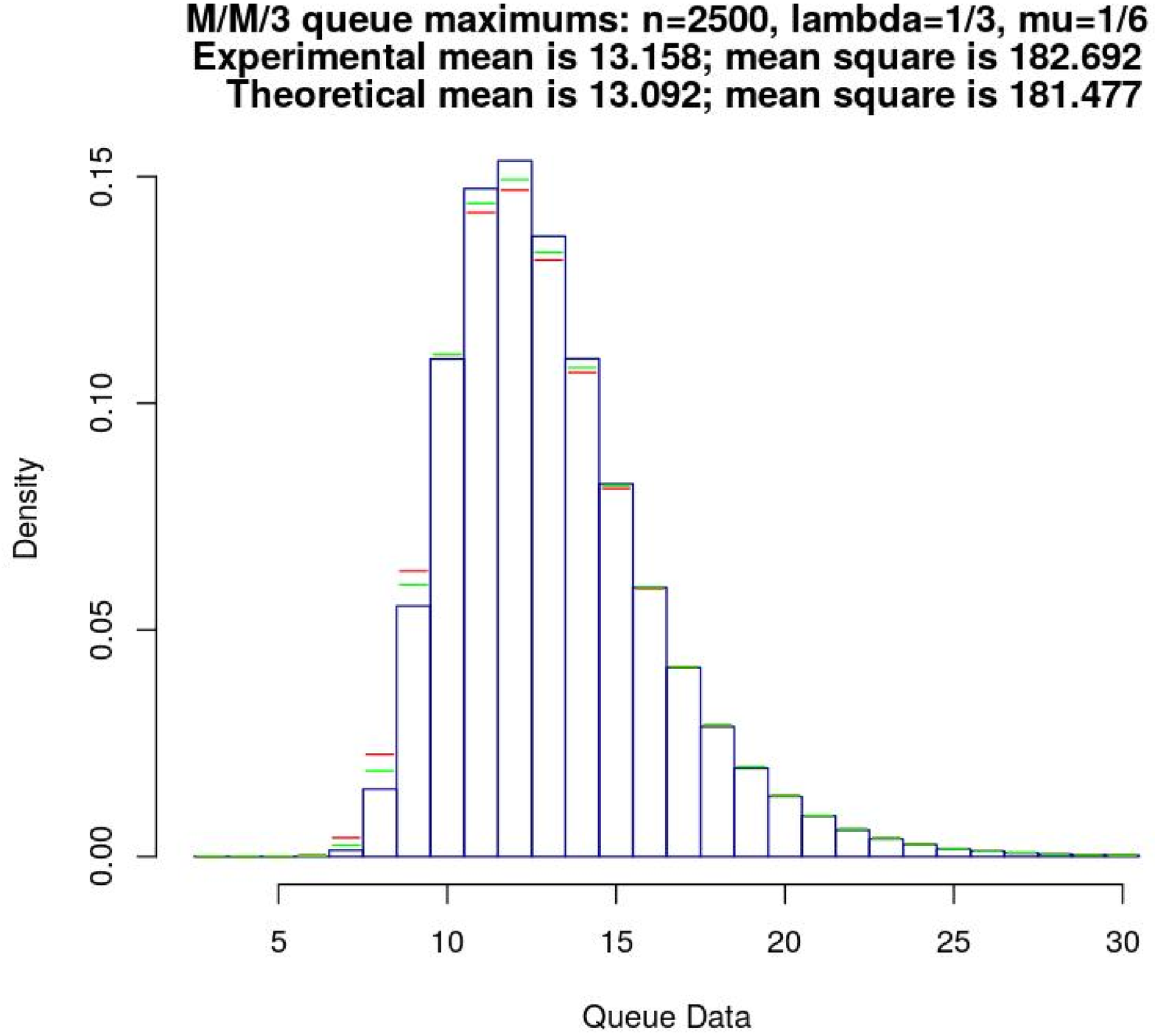}%
\end{figure}
\begin{figure}[ptb]%
\centering
\includegraphics[
height=3.0234in,
width=3.4541in
]%
{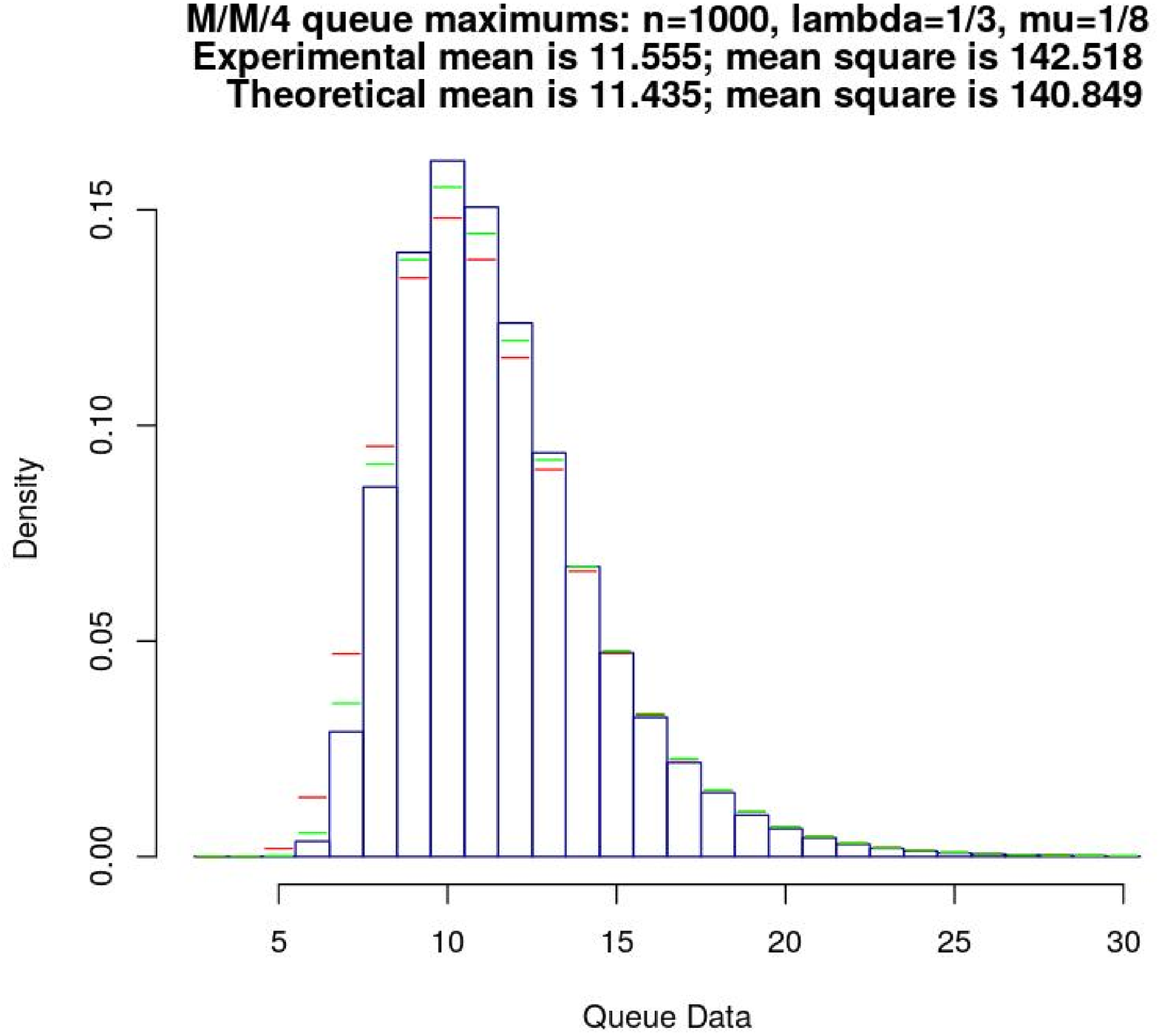}%
\end{figure}
\begin{figure}[ptb]%
\centering
\includegraphics[
height=3.0234in,
width=3.4541in
]%
{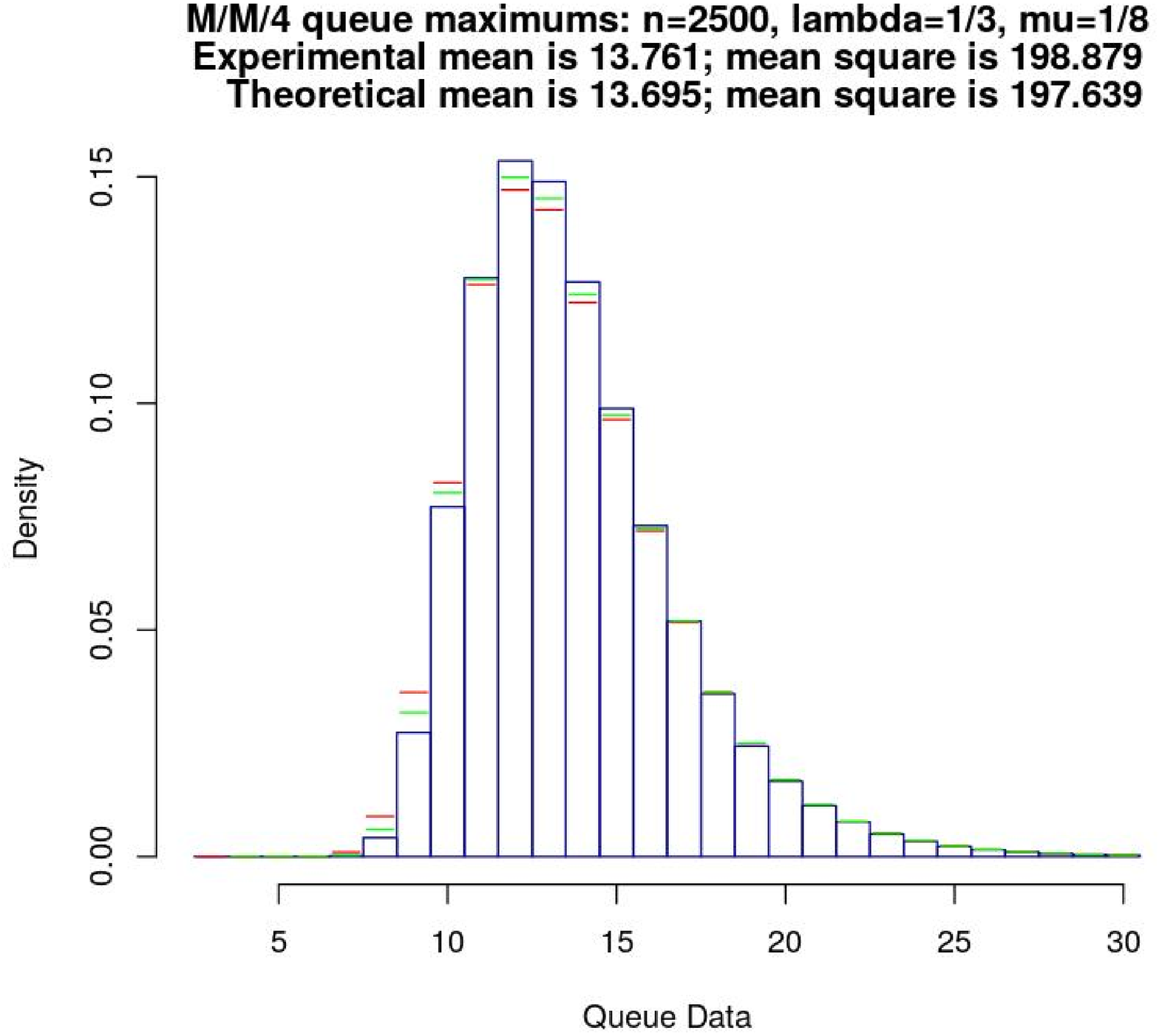}%
\end{figure}
\begin{figure}[ptb]%
\centering
\includegraphics[
height=3.0234in,
width=3.4541in
]%
{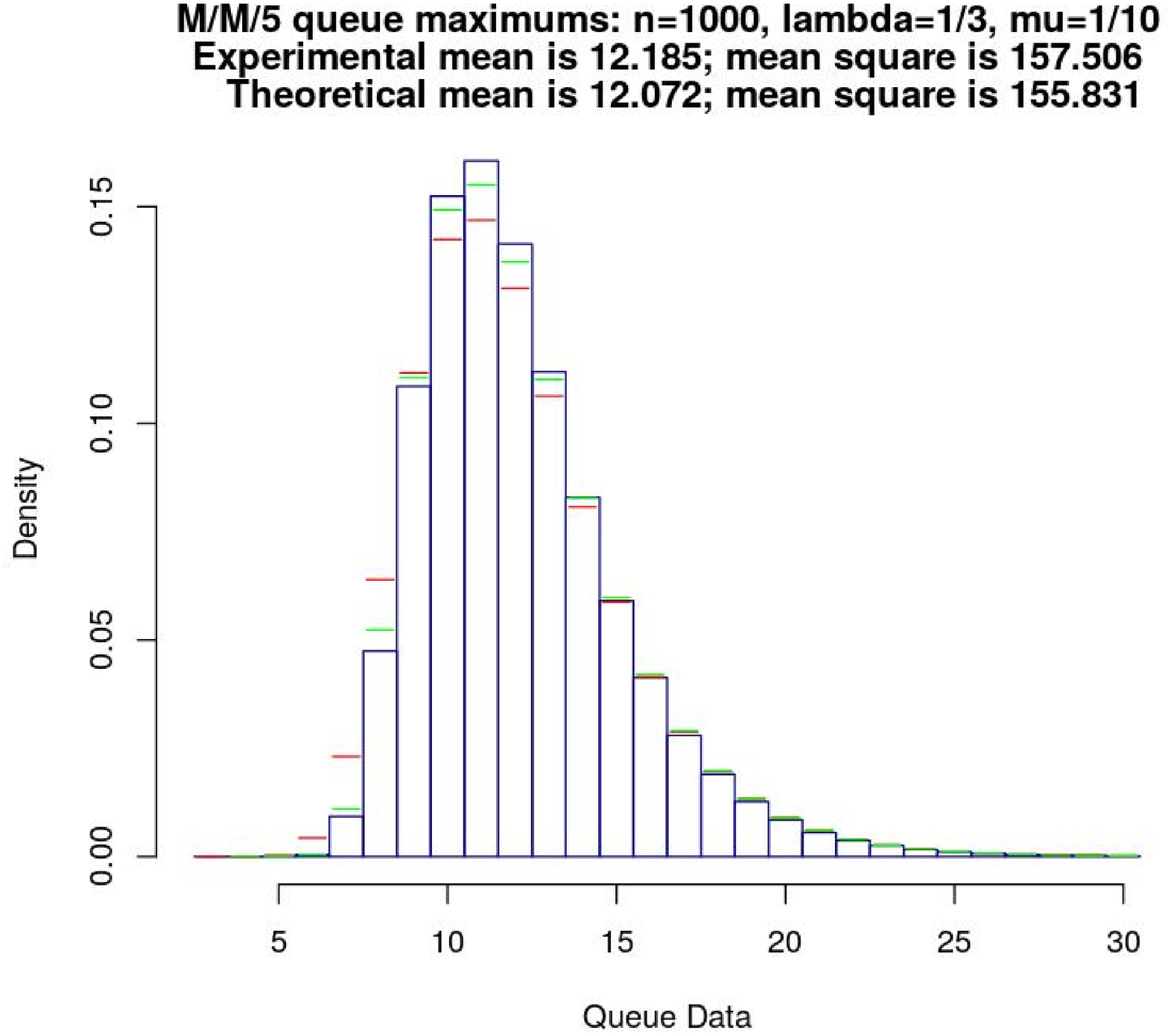}%
\end{figure}
\begin{figure}[ptb]%
\centering
\includegraphics[
height=3.0234in,
width=3.4541in
]%
{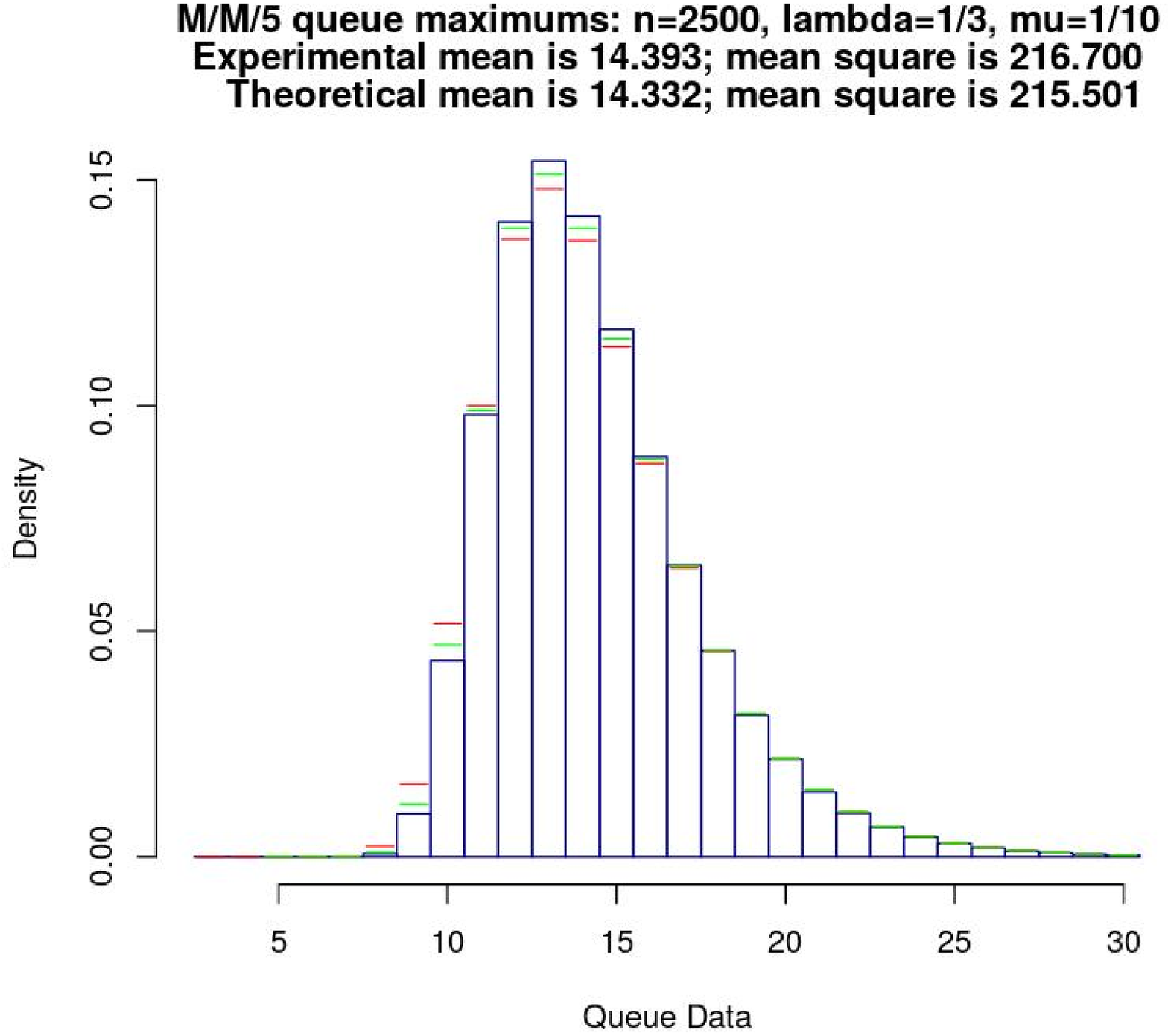}%
\end{figure}

\section{Erlang C}

Letting $\pi$ denote the stationary distribution of M/M/$c$, the probability
that all $c$ servers are busy is \cite{Cp-heu}%
\[%
{\displaystyle\sum\limits_{k=c}^{\infty}}
\pi_{k}=\frac{c^{c}}{c!}\left(
{\displaystyle\sum\limits_{k=c}^{\infty}}
\rho^{k}\right)  \pi_{0}=\frac{(c\rho)^{c}}{c!(1-\rho)}\pi_{0}
\]
where $\rho=\lambda/(c\mu)$ and
\[
\frac{1}{\pi_{0}}=%
{\displaystyle\sum\limits_{j=0}^{c-1}}
\frac{(c\rho)^{j}}{j!}+\frac{(c\rho)^{c}}{c!(1-\rho)}.
\]
We wish to demonstrate that%
\[
\frac{1}{\pi_{0}}=\frac{\sum\limits_{i=1}^{c}i!\tbinom{c-1}{i-1}\lambda
^{c-i}\mu^{i-1}}{(c-1)!\mu^{c-2}(c\mu-\lambda)}.
\]
As a preliminary step, note that%
\begin{align*}
\frac{(c\rho)^{c}}{c!(1-\rho)}  &  =\frac{\rho(c\rho)^{c-1}}{(c-1)!(1-\rho
)}=\frac{-(1-\rho)(c\rho)^{c-1}+(c\rho)^{c-1}}{(c-1)!(1-\rho)}\\
&  =-\frac{(c\rho)^{c-1}}{(c-1)!}+\frac{(c\rho)^{c-1}}{(c-1)!(1-\rho)}.
\end{align*}
The new $(1/\pi_{0})$-formula is equal to%
\[
\frac{1}{c\mu}%
{\displaystyle\sum\limits_{i=1}^{c}}
\frac{i}{(c-i)!}\frac{\lambda^{c-i}\mu^{i-c+1}}{1-\lambda/(c\mu)}=\frac{1}{c}%
{\displaystyle\sum\limits_{i=1}^{c}}
\frac{i}{(c-i)!}\frac{(c\rho)^{c-i}}{1-\rho}.
\]
As index $i$ runs from $1$ to $c$, index $j=c-i$ runs from $c-1$ to $0$,
giving%
\[
\frac{1}{c}%
{\displaystyle\sum\limits_{j=0}^{c-1}}
\frac{c-j}{j!}\frac{(c\rho)^{j}}{1-\rho}=%
{\displaystyle\sum\limits_{j=0}^{c-1}}
\frac{(c\rho)^{j}}{j!(1-\rho)}-\rho%
{\displaystyle\sum\limits_{j=1}^{c-1}}
\frac{(c\rho)^{j-1}}{(j-1)!(1-\rho)}
\]
which telescopes to%
\[%
{\displaystyle\sum\limits_{j=0}^{c-2}}
\frac{(c\rho)^{j}}{j!}+\frac{(c\rho)^{c-1}}{(c-1)!(1-\rho)}=%
{\displaystyle\sum\limits_{j=0}^{c-1}}
\frac{(c\rho)^{j}}{j!}-\frac{(c\rho)^{c-1}}{(c-1)!}+\frac{(c\rho)^{c-1}%
}{(c-1)!(1-\rho)}.
\]
By the preliminary step, this collapses to the old $(1/\pi_{0})$-formula and
we are done.

\section{Erlang B}

In an M/M/$c$/$c$ queue, if a customer arrives when all $c$ servers are busy,
the customer leaves the system immediately (with no effect on the queue).
\ The probability that all $c$ servers are busy is \cite{Cp-heu}%
\[
\pi_{c}=\frac{(c\rho)^{c}}{c!}\pi_{0}
\]
where $\rho=\lambda/(c\mu)$ and%
\[
\frac{1}{\pi_{0}}=%
{\displaystyle\sum\limits_{j=0}^{c}}
\frac{(c\rho)^{j}}{j!}.
\]

\section{Erlang A}

In an M/M/$c$+M queue, customers arrive with patience times $\tau$ that are
independent, exponentially distributed with mean $1/\theta$. The abandonment
rate $\theta$ is $0$ for Erlang C and is $\infty$ for Erlang B. \ If no
service is offered before time $\tau$ has elapsed, the customer leaves the
system immediately. \ Define%
\[
E=\frac{\frac{(c\rho)^{c}}{c!}}{\sum\limits_{j=0}^{c}\frac{(c\rho)^{j}}{j!}},
\]
the ratio from Section 2; and for $x>0$, $y\geq0$,%
\[
A(x,y)=\frac{x\exp(y)}{y^{x}}%
{\displaystyle\int\limits_{0}^{y}}
t^{x-1}\exp(-t)dt=1+%
{\displaystyle\sum\limits_{k=1}^{\infty}}
\frac{y^{k}}{\prod\limits_{\ell=1}^{k}(x+\ell)},
\]
an incomplete gamma function. \ The probability that all $c$ servers are busy
is \cite{P-heu, Z-heu, R-heu}%
\[%
{\displaystyle\sum\limits_{k=c}^{\infty}}
\pi_{k}=\frac{(c\rho)^{c}}{c!}A\left(  \frac{c\mu}{\theta},\frac{\lambda
}{\theta}\right)  \pi_{0}
\]
and \
\[
\frac{1}{\pi_{0}}=\frac{(c\rho)^{c}}{c!}\left[  \frac{1}{E}+A\left(
\frac{c\mu}{\theta},\frac{\lambda}{\theta}\right)  -1\right]  .
\]
We wonder about the implications of work in \cite{PW-heu}, especially a result
involving the constants $\gamma$ and $\pi^{2}/6$. \ Might certain issues we've
neglected here concerning asymptoptic moments ($h$ is real, not integer) be resolvable?

\section{Acknowledgements}

I am thankful to Guy Louchard for introducing me to the Poisson clumping
heuristic and to Stephan Wagner for extracting discrete Gumbel asymptotics in
\cite{Fc-heu} (a contribution leading to \cite{Fi1-heu, Fi0-heu} and the
present work). \ Writing simulation code for M/M/$c$ was facilitated by a
theorem in \cite{Sig-heu} involving order statistics of iid Uniform rvs. \ The
creators of R and Mathematica, as well as administrators of the
MIT\ Supercloud Cluster, earn my gratitude every day.


\begin{thebibliography}{99}                                                                                               %


\bibitem {Ald-heu}D. Aldous, \textit{Probability Approximations via the
Poisson Clumping Heuristic}, Springer-Verlag, 1989, pp. 1--8, 23--25, 30; MR0969362.

\bibitem {Fi1-heu}S. Finch, Geo/Geo/$2$ queues and the Poisson clumping
heuristic, arXiv:1902.09272.

\bibitem {Fn-heu}S. R. Finch, Euler-Mascheroni constant, \textit{Mathematical
Constants}, Cambridge Univ. Press, 2003, pp. 28--40; MR2003519.

\bibitem {Sfz-heu}R. F. Serfozo, Extreme values of birth and death processes
and queues, \textit{Stochastic Process. Appl.} 27 (1988) 291--306; MR0931033.

\bibitem {MP-heu}W. P. McCormick and Y. S. Park, Approximating the
distribution of the maximum queue length for M/M/$s$ queues, \textit{Queueing
and Related Models}, ed. U. Narayan Bhat and I. V. Basawa, Oxford Univ. Press,
1992, pp. 240--261; MR1210568.

\bibitem {HM-heu}G. Hooghiemstra, L. E. Meester and J. H\"{u}sler, On the
extremal index for the M/M/$s$ queue, \textit{Comm. Statist. Stochastic
Models} 14 (1998) 611--621; MR1621358.

\bibitem {Cp-heu}R. B. Cooper, \textit{Introduction to Queueing Theory},
2$^{\text{nd}}$ ed., North-Holland, 1981, pp. 79--101, 176--178; MR0636094.

\bibitem {P-heu}C. Palm, Research on telephone traffic carried by full
availability groups, \textit{Tele} 1 (1957) 1--107 [Engl. transl. of five
papers published in Swedish in 1946].

\bibitem {Z-heu}S. Zeltyn, \textit{Call Centers with Impatient Customers:
Exact Analysis and Many-Server Asymptotics of the M/M/}$n$\textit{+G Queue},
Ph.D. thesis, Israel Institute of Technology, 2004;
http://ie.technion.ac.il/\symbol{126}serveng/course2004/References/MMNG\_thesis.pdf.

\bibitem {R-heu}L. Rozenshmidt, \textit{On Priority Queues with Impatient
Customers: Stationary and Time-Varying Analysis}, M.Sc.\ thesis, Israel
Institute of Technology, 2007; http://ie.technion.ac.il/\symbol{126}serveng/course2004/References/thesis\_Luba\_Eng.pdf.

\bibitem {PW-heu}G. Pang and W. Whitt, Heavy-traffic extreme value limits for
Erlang delay models, \textit{Queueing Syst.} 63 (2009) 13--32; MR2576005;
http://www.columbia.edu/\symbol{126}ww2040/PangWhittExtremesErlang.pdf.

\bibitem {Fc-heu}S. Finch, The maximum of an asymmetric simple random walk
with reflection, arXiv:1808.01830.

\bibitem {Fi0-heu}S. Finch and G. Louchard, Traffic light queues and the
Poisson clumping heuristic, arXiv:1810.12058.

\bibitem {Sig-heu}K. Sigman, Notes on the Poisson process, lecture notes
(2009), http://www.columbia.edu/\symbol{126}ks20/stochastic-I/stochastic-I-PP.pdf.%

\begin{tabular}
[c]{lll}
& Steven Finch & \\
& MIT Sloan School of Management & \\
& Cambridge, MA, USA & \\
& \textit{steven\_finch@harvard.edu} &
\end{tabular}

\end{thebibliography}
\end{document}